\documentclass[conference]{IEEEtran}

\usepackage{hyperref}
\usepackage[numbers,sort]{natbib}
\usepackage{amsmath,amsfonts,amsthm, amssymb}
\usepackage{subfiles,xspace,cleveref,stmaryrd,microtype}
\usepackage{bm, derivative, multirow}

\usepackage{todonotes,subcaption}
\setuptodonotes{inline}

\newcommand{\filled}[1]{{\color{black} #1}}

\newcommand{\veps}{\varepsilon}

\newcommand{\acos}{\operatorname{acos}}
\newcommand{\asin}{\operatorname{asin}}
\newcommand{\atan}{\ensuremath{\operatorname{atan}}}
\newcommand{\pow}{\operatorname{pow}}
\newcommand{\fmod}{\ensuremath{\operatorname{fmod}}\xspace}
\newcommand{\logp}{\ensuremath{\operatorname{log1p}}\xspace}
\newcommand{\hypot}{\ensuremath{\operatorname{hypot}}\xspace}
\newcommand{\copysign}{\ensuremath{\operatorname{copysign}}\xspace}
\newcommand{\name}{\textsc{ExplaniFloat}\xspace}
\newcommand{\numbench}{546\xspace}

\newcommand{\highvarnum}{16\xspace}
\newcommand{\avgvarnum}{2.7\xspace}

\newcommand{\highopnum}{360\xspace}
\newcommand{\avgopnum}{9.5\xspace}
\newcommand{\halfprecnum}{33\xspace}
\newcommand{\noerrornum}{115\xspace}

\begin{document}

%% Title information
\title{Mixing Condition Numbers and Oracles \\ for Accurate Floating-point Debugging}

\author{
\IEEEauthorblockN{Bhargav Kulkarni}
\IEEEauthorblockA{Kahlert School of Computing\\
University of Utah \\
\texttt{bhargavk@cs.utah.edu}}
\and
\IEEEauthorblockN{Pavel Panchekha}
\IEEEauthorblockA{Kahlert School of Computing\\
University of Utah\\
\texttt{pavpan@cs.utah.edu}}}

\maketitle

\begin{abstract}
  Recent advances have made
    numeric debugging tools much faster
    by using double-double oracles,
    and numeric analysis tools much more accurate
    by using condition numbers.
  But these techniques have downsides:
    double-double oracles have correlated error
    so miss floating-point errors
    while condition numbers
    cannot cleanly handle over- and underflow.
  We combine both techniques to avoid these downsides.
  Our combination, \name,
    computes condition numbers using double-double arithmetic,
    which avoids correlated errors.
  To handle over- and underflow,
    it introduces a separate logarithmic oracle.
  As a result, \name achieves a precision of \filled{80.0\%}
    and a recall of \filled{96.1\%}
    on a collection of \numbench difficult numeric benchmarks:
    more accurate than double-double oracles
    yet dramatically faster than
    arbitrary-precision condition number computations.
\end{abstract}


\section{Introduction}
\label{sec:intro}

Floating-point numbers approximate real values
  but introduce subtle errors that can be difficult to detect,
  sometimes with catastrophic consequences.
Two prominent classes of automated tools
  attempt to improve the situation.
Numeric debugging tools~%
  \cite{fpdebug,herbgrind,fpsan,eftsan,shaman,verrou}
  observe numeric program executions
  and warn the programmer about error-inducing operations.
By contrast, static analysis tools~%
  \cite{salsa,daisy,fptaylor,satire,rosa,fluctuat,gappa}
  analyze a short numeric program's behavior
  over an entire region of possible inputs
  and provide varyingly-sound guarantees
  about the maximum possible floating-point error.

Recent years have seen rapid improvement in both tool classes.
In debugging, a line of tools from Herbgrind~\cite{herbgrind}
  to FPSanitizer~\cite{fpsan} and EFTSanitizer~\cite{eftsan}
  has focused on reducing runtime overhead
  by introducing more optimized ways to compute
  \textit{oracle values} for floating-point values.
Herbgrind computes these oracles using a JIT-compiled virtual machine,
  FPSanitizer uses a compiler pass to
  insert native calls to the GNU~MPFR library,
  and EFTSanitizer replaces MPFR
  with inlined double-double computations.
These innovations reduce overhead from
  $574\times$ to $111\times$ to $12.3\times$,
  but also reduce accuracy and cause false negatives.
In static analysis, a line of tools from Salsa~\cite{salsa}
  to Rosa~\cite{rosa} and FPTaylor~\cite{fptaylor}
  have focused on improving error estimation
  through more-accurate representations of \textit{error bounds}.
Salsa uses value and error interval arithmetic,
  Rosa uses affine arithmetic,
  and FPTaylor uses error Taylor series.
These innovations have dramatically tightened
  achievable error bounds,
  but the best techniques cannot reason accurately
  about overflow (and, for some tools, underflow).
More importantly,
  the performance innovations of debugging tools
  and accuracy innovations of static analysis tools
  have not yet been combined to achieve
  both performance and accuracy simultaneously.

% Here's how our solution works
This paper introduces \name,%
\footnote{
\name is open source and available online at
\url{https://github.com/herbie-fp/herbie/tree/bhargav-nobigfloat}.
}
  a floating-point debugging tool that combines
  double-double oracle values~\cite{eftsan}.
  with condition number Taylor error bounds~\cite{fptaylor,atomu}
  to detect erroneous operations.
It also accurately detects over- and underflow errors
  using a logarithmic oracle for out-of-range values.
Our implementation based on the \texttt{qd} library~\cite{qd}
  achieves high accuracy and high performance
  using a novel implementation of these oracles:
  on \numbench benchmarks from the Herbie~2.1 suite,
  \name achieves a a precision of \filled{80.0\%}
  and a recall of \filled{96.1\%}.
By contrast, a traditional double-double oracle achieves
  a precision of \filled{56.5\%} and recall of \filled{65.4\%}.
\name's accuracy is comparable to an arbitrary-precision baseline,
  while also being significantly ($4.24\times$) faster.
In short, this paper contributes:

\begin{itemize}
\item A debugging algorithm based on condition numbers instead of
  oracle values (\Cref{sec:condition}).
\item A novel logarithmic oracle for accurately tracking overflow and underflows (\Cref{sec:overflow}).
\item An implementation using a novel number representation and the \texttt{qd} library (\Cref{sec:dsl}).
\end{itemize}

\section{Background and Related Work}
\label{sec:background}

As is well known, floating-point numbers $\hat x$ are
  a subset of the real numbers $x$
  with a fixed precision and exponent range.
We write $R(x)$ for
  the closest floating-point number to the real number $x$;
  $R$ suffers from \emph{rounding error} due to limited precision
  and \emph{over- and underflow} due to limited exponent range.

\subsection{Debugging, Oracles, and Shadow Memory}
\label{subsec:oracle}

A \emph{numerical error} occurs
  when the floating-point result of a program
  differs from the correct real-number result.
Modern debugging tools detect these errors using \emph{oracles}:
  for every floating-point intermediate $\hat{x}_i$
  they store an \emph{oracle value} $x_i$ in higher precision.
By comparing $\hat{x}_i$ to $x_i$ they can then detect numerical errors.
Examples of this approach include FpDebug~\cite{fpdebug},
  Herbgrind~\cite{herbgrind}, and FPSanitizer~\cite{fpsan},
  while Shaman~\cite{shaman} and Verrou~\cite{verrou}
  use a similar technique but probabilistically.
The state of the art is the recent EFTSanitizer tool,
  which uses an oracle based on
  \emph{double-double arithmetic}.

A double-double value is a pair $(\hat{x}, \hat r)$
  representing the real value $x \approx \hat x + \hat r$.
Applying a function $\hat{f}_i(\hat{x}_j)$ in double-double
  yields $\hat{x}_i$ and $\hat{r}_i$ such that
\begin{align*}
  \hat{x}_i &= \hat{f}_i(\hat{x}_j) \\
  \hat{r}_i &= R(f_i(\hat{x}_j + \hat{r}_j) - \hat{f}_i(\hat{x}_j))
\end{align*}
The intuition is that $\hat{r}_i$ captures the error of $\hat{x}_i$,
  providing, in effect, extra precision for $\hat{x}_i$.
Importantly, $\hat{r}_i$ can be computed
  using only hardware floating-point instructions
  much faster than an arbitrary-precision libraries like MPFR.

However, double-double oracles
  can fall prey to the same rounding errors
  as the original computation.
For example, in the textbook example~\cite{book87-nmse}
  $\sqrt{x+1} - \sqrt{x}$, thousands of bits of precision are needed
  to accurately compute the result~\cite{herbgrind,herbie}.
The double-double oracle doesn't have enough precision
  and thus computes the same erroneous value
  as the original computation.
This causes a numerical debugger using a double-double oralce
  to miss the numerical error---a false negative.
Double-double oracles also cannot handle overflow and underflow,
  since double-double values have an exponent range
  no larger than ordinary double-precision floating-point.

\subsection{Error Taylor Series and Condition Numbers}

Numerical errors can be quantified via the relative error:
  a program $\hat P$ has a \emph{relative error bound} $c$
  if $|\hat{P}(x) - P(x)|$ is bounded by $c| P(x)|$
  for all $x$ in some set of inputs.%
\footnote{
  The case when $P(x) = 0$ requires separate,
  usually tool-specific, handling.
}
Modern error analysis tools
  derive such relative error bounds by
  composing known error bounds
  $\hat{f}(x) = f(x)(1 + \varepsilon)$ (where $|\varepsilon| < c$)
  for primitive operations $f(x)$..
Examples of such tools include
  Salsa~\cite{salsa}, Rosa~\cite{rosa}, Daisy~\cite{daisy},
  Fluctuat~\cite{fluctuat}, Gappa~\cite{gappa},
  Precisa~\cite{precisa}, and FPTaylor~\cite{fptaylor}.
Note that while debugging tools consider a single input at a time,
  worst-case error bound tools aim to reason
  abstractly about a set of possible inputs.
The state of the art is the recent Satire~\cite{satire} tool,
  which uses \emph{error Taylor series} computed using
  \emph{automated differentiation}.%
\footnote{Satire's approach is different from,
  but closely related to condition numbers.}

Error Taylor series replace every primitive operation $f(x)$
  in the program $\hat{P}(x)$ by $f(x)(1 + \varepsilon)$,%
\footnote{Using a unique $\varepsilon$ for each operation.}
  resulting in a real-number formula $P(x, \varepsilon)$
  where each $\varepsilon$ has bounded magnitude.
Note that $P(x, 0) = P(x)$, the ideal real-number behavior;
  thus, a worst-case error bound for $\hat{P}$ can be derived
  by studying how $P(x, \varepsilon)$ varies in $\varepsilon$.
We can estimate the worst-case error by taking a Taylor expansion:
\begin{multline*}
	\hat P(x, \veps) =
        \underbrace{P(x, 0)}_{\text{Exact answer}} +
        \underbrace{ \sum_{\veps} \veps \left. \pdv{P}{\veps} \right|_{x}}_{\text{First-order error}} + \\
        \underbrace{ \sum_{\veps}\sum_{\veps'}
        \veps\veps' \left. \pdv{P}{\veps, \veps'} \right|_{x} + \dotsb}_{\text{Higher-order error}}
\end{multline*}
The second-order error is typically small~\cite{fptaylor}
  and usually ignored~\cite{satire};
  subtracting off the exact answer then leaves:
\begin{equation}
  \label{eq:foe}
c = \left|\sum_i
          \veps_i
          \left. \pdv{P}{\veps_i} \right|_{x}
          \frac{1}{P(x)}
          \right|
     <
    \sum_i c_i
          \underbrace{
    \left|
          \left. \pdv{P}{\veps_i} \right|_{x}
          \frac{1}{P(x)}
    \right|}_{A_i(x)}
\end{equation}
  where $c_i$ is the maximum relative error
  of the $i$-th primitive operation in $\hat{P}$.
Each $A_i$ is a real-valued function of $x$
  and can be bounded using interval arithmetic,
  global non-linear optimization, or other techniques
  to compute a worst-case error bound for $\hat P$.
The important take-away here
  is that worst-case error bounds can be computed
  \emph{without} computing the exact value $P(x)$
  and thus without running the risk that that exact value
  will be computed incorrectly,

Condition numbers are a convenient shortcut
  for computing this \Cref{eq:foe}.
The condition number $\Gamma_f(x)$ of a computation $f(x)$ is
\[
  \Gamma_f(x) = \left| \frac{x f'(x)}{f(x)} \right| = \left|\frac{\partial f / f}{\partial x / x} \right|
\]
The condition number has the following property:
  if $\hat x = x (1 + \veps)$,
  then $f(\hat x) = f(x)(1 + \Gamma_f(x)\veps) + O(\veps^2)$,
  where the $O(\veps^2)$ term can be ignored
  when computing only the first-order error.
In other words, $\Gamma_f$ measures
  how much $f$ amplifies incoming relative error.
Note that the condition number is defined
  purely in terms of the real-number function $f$;
  it is an inherent property of the function
  conserved across precisions
  and not dependent on the quality of an oracle.
The relative error of each intermediate value
  $x_i = f_i(x_j)$ in $\hat P$
  can then be computed with
\begin{align*}
x_i &= f_i(x_j) \\
|E_i| &= \Gamma_{f_i}(x_j) |E_j| + c_i
\end{align*}
  and an analogous formula for binary functions.
We refer the reader to the ATOMU paper~\cite{atomu}
  for a more readable and detailed derivation;
  the main take-away is that the first-order method
  can be implemented by simply computing,
  for each intermediate value $x_i$,
  an error bound $E_i$.
In practice, this produces accurate error bounds
  even when the intermediate value $x_i$
  suffers from rounding error.

\section{Condition Numbers for Rounding Error}
\label{sec:condition}

\name is based on the observation
  that modern debugging and error-bound approaches
  both execute a program
  while also tracking additional metadata ($|E_i|$ and $\hat{r}$)
  that estimates the error of that computation.
It thus uses a hybrid approach,
  executing the floating-point program,
  but using condition numbers to detect rounding error.
Specifically, \name executes the floating-point program
  and computes a double-double oracle
  for each floating-point intermediate.
However, it detects possible numerical errors
  not by comparing the actual and oracle value
  but by computing the condition number of each operation.

This provides two advantages.
Firstly, since the condition number does not depend on an ``oracle'',
  it is robust to innaccuracies in the oracle.
In fact, in our evaluation (\Cref{sec:eval}),
  \name achieves much better precision and recall
  than a similar tool using the oracle method,
  largely due to innacuracies of the oracle.
Secondly, though harder to evaluate,
  condition numbers have better error localization.
Each condition number is computed from
  a specific floating-point operation,
  and every warning raised by \name indicates this operation.
Comparing the oracular and computed values,
  by contrast, implicates the full program execution up to that point
  so may raise errors too late, or warn the user
  about the slow and steady accumulation of error
  over a long series of slightly-erroneous computations
  where any single operation is a red herring.

In practice, we found that condition numbers on their own
  make for a poor debugging experience.
A function like $\sin$ has high condition numbers
  both for large inputs and for inputs close to a multiple of $\pi$;
  merely identifying the problematic operation
  didn't give users all the information they needed.
\name therefore \emph{splits} the condition numbers
  of each supported operation.
For example, the standard condition number for $\sin(x)$
  is $\Gamma_{\sin}(x) = |x \cot(x)|$; we split this into two parts:
  $\Gamma^1_{\sin}(x) = |x|$
  and
  $\Gamma^2_{\sin}(x) = |\cot(x)|$,
  with $\Gamma_{\sin} = \Gamma^1_{\sin} \Gamma^2_{\sin}$.
$\Gamma^1_{\sin}(x)$ indicates ``stability'' errors for $\sin(x)$,
  while $\Gamma^2_{\sin}(x)$ indicates ``cancellation'' errors,
  where ``stability'' refers to errors for very large or small inputs
  while ``cancellation'' refers to errors for inputs
  close to some discrete set.
The split condition numbers for each operation
  supported by \name
  are shown in \Cref{tbl:conditions}.

\begin{table}[tbp]
\centering
\begin{tabular}{l | lll}
Operation &
Condition Number &
Bad inputs &
Type \\\hline

$x \pm y$ &
$|\{x, y\} / (x \pm y)|$ &
$x \approx y$ &
Cancellation \\

$x \cdot y$, $x / y$ &
- &
- &
- \\\hline

$\sqrt{x}$, $\sqrt[3]{x}$ &
$\frac12$, $\frac13$ &
- &
- \\\hline

$\log(x)$ &
$| 1 / \log(x) |$ &
$x \approx 1$ &
Cancellation \\

$\exp(x)$ &
$|x|$ &
$x$ large &
Sensitivity \\

$x^y$ &
$|y|$ and $| y \log x |$ &
$y$ large &
Sensitivity \\\hline

$\sin(x)$ &
$|1 / \tan(x)|$ &
$x \approx k \pi$ &
Cancellation \\

&
$|x|$ &
$x$ large &
Sensitivity \\

$\cos(x)$ &
$|\tan(x)|$ &
$x \approx (k + \frac12) \pi$ &
Cancellation \\

&
$|x|$ &
$x$ large &
Sensitivity \\

$\tan(x)$ &
$|\tan(x) + 1/\tan(x)|$ &
$x \approx (\frac{k}2) \pi$ &
Cancellation \\

&
$|x|$ &
$x$ large &
Sensitivity \\\hline

$\acos(x)$ &
$|x / (\sqrt{1 - x^2}\acos(x)) |$        &
$|x| \approx 1$ &
Cancellation \\

$\asin(x)$ &
$| x / (\sqrt{1 - x^2}\asin(x)) |$ &
$|x| \approx 1$ &
Cancellation \\
\end{tabular}

\caption{
  All operations and split condition numbers supported by \name.
  Dashes indicate unused error types
    for a particular operation.
  Note that all cancellation errors
    are caused by the input $x$ being close to
    some specific value (or discrete set of values)
    while all sensitivity errors
    are caused by large inputs.
}
\label{tbl:conditions}
\end{table}

All told, \name detects and warns the user about rounding error
  for each operation with high condition number,
  testing each split condition number for that operation
  and warning the user for each one that is past
  a user-specified threshold.
However, there is an exception to this rule
  in the case of operations on exact values, like in $\log(1)$.
Here, the condition number $|1 / \log(1)|$ is infinite,
  meaning that the $\log$ operation significantly amplifies any input error;
  however, because $1$ is exactly-represented in floating-point,
  it has no error to amplify.
To avoid raising false alarms in such situations,
  \name does not warn for high condition numbers
  for arguments that are exact constants.
For multi-argument functions where some arguments are exact,
  only condition numbers associated with non-exact arguments
  produce a warning.

Despite these tweaks, there is still one common cause
  of false postives: operations that introduce minimal error.
Normally, floating-point operation introduce
  around one machine epsilon of error.
Some operations, however, introduce much less.
For example, the expression $2^{100} + 2^{-100}$
  evaluates to $2^{100}$, with $2^{-200}$ relative error.
Later computations can amplify that error by a lot
  and still not have significant error;
  for example, in $\cos(2^{100} + 2^{-100})$
  the condition number $2^{100}$
  only amplifies the error to $2^{-200} 2^{100} = 2^{-100}$.
This issue occurs rarely (see \Cref{sec:eval}),
  but is a notable case where comparing
  actual and oracle values would be more accurate.

\section{Oracles for Overflow and Underflows}
\label{sec:overflow}

While condition numbers work better than oracles
  for rounding error, we found the opposite to hold
  for over- and underflow errors.
In fact, both the debugging and static analysis literature
  treat any overflow or underflow as an error,%
\footnote{
FPTaylor~\cite{fptaylor}
  does have specialized handling for underflow
  (via its $f(x)(1 + \varepsilon) + \delta$ error model)
  but treats any overflowing operation as an error.
}
  an approach reminiscent of treating all large condition numbers as errors.
This weak modeling of over- and underflows
  causes false positives, for example,
  in $1 + 1/\exp(x)$, $\exp(x)$ overflows for large $x$
  but the full expression still correctly evaluates to $1$.

Condition numbers cannot help with over- and underflows:
  condition numbers are based on relative error bounds for primitive operations
  that do not hold when overflow (and sometimes underflow) occurs.
But a oracle \emph{could} avoid false positives:
  instead of raising a warning when an expression overflows,
  an over- and underflow oracle
  would approximate the overflowed value
  and track whether the overflow actually caused the computation
  to diverge from a real execution.

Since an over- and underflow oracle requires
  a vastly larger dynamic range than ordinary floating-point,
  \name uses a logarithmic number system as an oracle.
In this system, a real number $x$ is represented
  by its sign plus the floating-point number $R(\log_2(|x|))$.
Even extremely large numbers are representable directly.%
\footnote{
The logarithmic number system can itself overflows,
  but this doesn't happen in our evaluation suite.
}
Operations $f(x)$ on oracle values
  require computing $\log_2(|f(2^x)|)$;
  for example, to compute an oracle for $\pow(x, y)$ (for positive $x$)
  one instead computes $y\log_2(x)$.
Luckily, such an oracle can be implemented efficiently using
  only hardware floating-point operations (see \Cref{sec:dsl}).

\name uses this logarithmic representation
  to approximate values outside the standard floating-point range,
  and warns when these out-of-range values cause
  the real and floating-point computation to diverge.
Consider the expression $\sqrt{1 + x^2}$ in double precision;
  the $x^2$ term can overflow for very large $x$ like $10^{300}$..
In this case the exact real-number value of $x^2$ is $10^{600}$
  while the floating-point result is $+\infty$.
But, since $+\infty$ is in fact the best double-precision representation
  of $10^{600}$, the floating-point and real executions
  have not yet actually diverged and no error is raised.
Instead, \name just represents the value
  logarithmically as $\log_2(10^{600}) \approx 1993$.
However, $\sqrt{1 + x^2}$ evaluates to $+\infty$ in floating-point
  while the logarithmic oracle is approximately $10^{300}$.
These \emph{are} starkly different in double-precision,
  so at this point \name determines that the
  floating-point and real executions have diverged
  and raises an error.

In general \name raises an error for any
  operation on over- or underflowed values
  whose oracle result is within the standard floating-point range.
One particular class of underflow errors, however
  is an exception from this rule.
For example, consider the same expression $\sqrt{1 + x^2}$
  but now for a very small value like $x = 10^{-300}$.
The addition operation has an out-of-range input ($x^2$ underflows)
  and produces an output in the standard range ($1$),
  but the underflow is benign because the true value, $1 + 10^{-600}$,
  still rounds to $1$.
To avoid such false positives,
  \name special-cases additions and subtractions where the two arguments
  whose differ significantly in order of magnitude,
  and over- and underflow errors in the smaller value are ignored.

\name's over- and underflow oracle again shows
  that combining both oracle and analytic techniques
  can reduce false positives
  without using arbitrary-precision floating-point.

\section{Implementation}
\label{sec:dsl}

The \name implementation aims to test
  the idea of a first-order method debugger
  with an oracle for overflow and underflow detection.
\name is thus aimed at debuggability rather than performance,
  and thus uses a simple floating-point virtual machine
  instead of the more complex (and performant)
  techniques pioneered in FPSanitizer and
  EFTSanitizer~\cite{fpsan,eftsan}.

\subsection{Shadow Memory}

Each floating-point intermediate in \name
  is shadowed by a double-double value,
  used for computing condition numbers,
  and a logarithmic value (plus a sign bit)
  for the over/underflow oracle:
\[
  \llbracket x \rrbracket = (\hat x, \hat r, s, \hat e)
\]
  where $\hat{x}$ is the floating-point value of $x$,
  $\hat{r}$ is the residual error of $x$,
  $s$ the sign of $x$,
  and $\hat e$ the logarithm of $|x|$.
We call this number system ``DSL'',
  after its components: double-double, sign, and logarithm.
The double-double shadow value
  improves the accuracy of computed condition numbers
  (reducing false positives) and also improves
  the accuracy of the over/underflow oracle
  after overflow or underflow occurs.

Every floating-point operation in our virtual machine
  performs shadow operations
  to compute the relevant $\hat x$, $\hat r$, $s$, and $\hat e$.
For values in the floating-point range,
  the shadow output's $\hat x$ and $\hat r$
  are computed from the shadow input's $\hat x$ and $\hat r$,
  and the $s$ and $\hat e$ values are computed from
  the shadow output.
For values outside the floating-point range,
  the input $\hat x$ and $\hat r$ values are ignored
  (since they typically contain zeros, infinities, and NaNs)
  and instead the shadow output's $s$ and $\hat e$
  are computed directly from the shadow input's $s$ and $\hat e$
  using standard logarithmic number system techniques.

After computing the shadow value,
  \name also computes the condition number of each operation
  using the same number representation.
If the condition number is greater than a user-configurable threshold,
  \name computes each split condition number
  and raises the appropriate error for each one over the threshold.
(Note that both the full and split condition numbers must be computed
  because functions like $\sin x$ can have split condition numbers
  ($x$ and $1 / \tan(x)$) that cancel out for some inputs ($x \approx 0$).)
Finally, if the output shadow value's $\hat e$ indicates that
  it is inside the representable floating-point range,
  while at least one input's shadow value is outside that range,
  \name also raises the appropriate over/underflow error,
  except in the case of suppression
  as described in \Cref{sec:condition,sec:overflow}.

Notably, all of these values are machine floating-point numbers
  (save $s$, which is a machine boolean)
  and do not require allocation;
  all shadow operations likewise use machine floating-point operations
  instead of arbitrary-precision arithmetic libraries
  like MPFR~\cite{mpfr}.
We expect this to enable high performance,
  much like in EFTSanitizer~\cite{eftsan}.
In fact, our evaluation (see \Cref{sec:eval})
  shows a $4.24\times$ speedup over an arbitrary precision library,
  and we expect much larger speedups in
  a real-world implementation focused on performance
  and run on larger benchmark programs.

The challenge with combining two number systems---%
  the double-double values for computing condition numbers
  and the logarithmic values for detecting problematic
  overflows and underflows---%
  is understanding how they interact.
In \name, only one component is active at a time:
  either the value is in the representable floating-point range,
  in which case the double-double value is used,
  or the value is out of the representable range
  and the logarithmic component is used;
  the logarithmic component determines which case a value is in.
As a result, a shadow operation may either
  receive inputs where the double-double component is active
  and produce logarithmic outputs, or vice-versa;
  both these cases require conversion.
Double-double inputs that overflow
  are first converted to logarithmic values,
  and then a logarithmic operation is performed;
  the conversion just requires taking the logarithm
  and recording the sign.
Logarithmic values that bring the value back in range
  compute an output double-double value
  by exponentiating.
Naturally, the resulting value is not precisely known,
  since the logarithmic component has less precision
  than even a simple double value, let alone a double-double one.
However,
  this logarithmic-to-double-double conversion
  only occurs after an over/underflow renormalization---%
  in other words, after an error is already detected---%
  so accuracy is less important.

\subsection{Implementing Double-Double Computation}

To implement double-double floating-point numbers,
  we wrap David Bailey's \texttt{qd} floating-point library~\cite{qd}.
This library provides double-double implementation \emph{cores}
  for a variety of operations.
These cores compute highly accurate $\hat x$ and $\hat r$ values,
  but typically only on a narrow range of values.
For example, the $\sin$ implementation in \texttt{qd}
  uses a naive range-reduction algorithm
  and is not accurate for inputs much larger than $\pi$.
The cores are also not robust to special values such as NaN;
  some cores like \atan\xspace even cause segmentation faults
  if called with NaN.
To address this,
  \name wraps \texttt{qd} with code that tests for special values
  and performs range reduction.

While handling special values was straightforward,
  range reduction took particular care,
  and our approach was heavily influenced by,
  and with significant code borrowed from,
  the classic \texttt{fdlibm} library~\cite{fdlibm}.
Consider the analog of $sin(x)$ in \texttt{qd}, \texttt{c\_dd\_sin},
  which is only accurate for inputs $x$ close to 0.
We need a wrapper around \texttt{c\_dd\_sin}
  to reduce input values to bring them closer to 0.
Since $\sin(x)$ is periodic, it is enough to subtract
  the relevant multiple of $2\pi$;
  in fact, it is enough to subtract the relevant multiple
  of $\pi/2$, and then dispatch to $\pm\sin(x)$ or $\pm\cos(x)$
  to compute the final value.
What's tricky is subtracting the relevant muliple of $\pi/2$
  in high precision from both the primary and residual parts
  of a double-double value.

To do so, we use the standard
  Payne-Hanek algorithm~\cite{payne-hanek}, as implemented in
  the \texttt{rem\_pio2} function from \texttt{fdlibm}.%
\footnote{
  For smaller input, \texttt{rem\_pio2} will use the
  faster Cody-Waite~\cite{cody-waite} algorithm.
}
This helper function converts an input $x$
  into outputs $x_1$, $x_2$, and $k$
  where $x \approx x_1 + x_2 + \frac{\pi}2 k$.%
\footnote{With $x_1$ and $x_2$ forming a
  double-double value and $k$ and integer between 0 and 4.}
To extend this approach to the double-double value $x + r$,
  we apply \texttt{rem\_pio2} to both $x$ and $r$,
  yielding
\[
x + r = \left(x_1 + x_2 + \frac\pi2 k_x\right) + \left(r_1 + r_2 + \frac\pi2 k_r\right)
\]
Then, $x_1 + x_2$ and $r_1 + r_2$ can be added as double-double values,
  yielding $x' + r'$; in other words,
\[
x + r = x' + r' + \frac\pi2 (k_x + k_r \bmod 4)
\]
The combined $x'$, $r'$, and $k_x + k_r \bmod 4$ values
  can then be passed to \texttt{c\_dd\_sin} (or \texttt{c\_dd\_cos})
  to compute $\sin(x + r)$.
In other words, we combine \texttt{fdlibm}'s range reduction
  with \texttt{qd}'s trigonometric cores
  to achieve a full-range double-double $\sin$ implementation.
A similar approach is used
  for other trigonometric functions.
For the logarithmic function,
  where a high-precision value of $e$ is required,
  we instead transplant a simpler implementation
  from the Racket \texttt{math/flonum} library's~\cite{racketmath}
  \texttt{fl2log} function.

The ultimate result is an elementary function library
  for double-double values that implements
  the shadow operations in \name.
Importantly, our implementation
  uses basically-standard double-precision algorithms
  for range reduction,
  and thus does not carry the performance penalty of
  an arbitrary-precision library.
To our knowledge, \name is
  the first numeric debugging tool
  with support for double-double transcendental
  function implementations.

\subsection{Implementing Logarithmic Computation}

To implement logarithmic floating-point numbers,
  we use the technique of~\citet{lns};
  in short, this technique centers around the use
  of $\Phi^+(x) = \log(1 + \exp(x))$
  and $\Phi^-(x) = \log(1 - \exp(x))$ functions
  for addition and subtraction.
Our implementation of these $\Phi$ functions
  internally uses the \texttt{qd} library
  for higher precision, in order to ensure accuracy.
Exponential, root, power, logarithm, and exponent operations
  use the straightforward implementation.
Some operations, such as trigonometric functions,
  are very challenging to implement
  for logarithmic values, so we don't try.
Instead, we note that applying these functions
  on out-of-range floating-point numbers
  returns in-range results (save near 0),
  meaning any use of these functions
  on over/underflowed inputs will raise
  an over/underflow error.
An accurate implementation is thus not needed.
Like the double-double operations,
  logarithmic operations use ordinary floating point operations,
  avoiding the overhead of arbitrary-precision computation.

\section{Evaluation}
\label{sec:eval}

We evaluate \name both in isolation as a numeric debugger
  and also by testing various components
  in an ablation study.
We focus on three research questions:

\begin{enumerate}
\item[\textbf{RQ1}]
  Does \name detect erroneous operations
    with few false positives and negatives?
\item[\textbf{RQ2}]
  Is \name more accurate than an oracle-based debugging tool?
\item[\textbf{RQ3}]
  Is \name as accurate as, but more performant than, an arbitrary-precision baseline?
\end{enumerate}

We chose the Herbie~2.1~\cite{herbie} benchmark suite as our evaluation target.
These 546 benchmarks%
\footnote{
Six benchmarks that use operations like \fmod, \logp, \hypot, and \copysign are not included in the results,
  as these operations are not supported by \name.
}
  are drawn from textbooks, papers, and open-source code
  and intended for evaluating floating-point repair tools,
  so have many complicated numerical errors,
  with many overflow, cancellation, and stability errors,
  many involving transcendental functions.
At the same time, they are small:
  up to \filled{\highvarnum} variables (\filled{\avgvarnum} average)
  and up to \filled{\highopnum} floating-point operations (\filled{\avgopnum} average) each.
This is critical for our evaluation,
  which evaluates \name by comparing to a baseline that uses
  very high arbitrary-precision (up to 10,000 bits)
  interval arithmetic to compute a ground truth.
These results should thus be indicative of
  the accuracy of \name's error detection;
  we expect \name's design,
  using only machine floating-point operations,
  to scale to much larger programs,
  but leave that evaluation for future work.

The Herbie~2.1 benchmarks come with a test runner
  that randomly samples 256 valid inputs for each benchmark;
  \filled{\noerrornum} benchmarks have no detected floating-point error,
  while the rest have at least some for some inputs.
Most benchmarks use 64-bit floating-point but some (\halfprecnum)
  use 32-bit floating-point; in either case \name uses
  64-bit floating-point for its shadow operations.
All experiments are run on a machine with a
  i7-8700K CPU (at 3.70GHz) and 32GB of DDR3 memory running
  Ubuntu~24.01, Racket~8.10, and MPFR version~4.2.1.

\subsection{RQ1: Predicting Floating-Point Error}

To determine \name's false positive and false negative rate,
  we need an accurate ground truth to compare to.
We compute one using the Rival interval arithmetic package~\cite{rival2}
  with up to 10\thinspace000-bit-precision floating-point.
\footnote{
  Five benchmarks are discarded
    because correctly-rounded evaluation fails.
}
If \name raises condition number or over/underflow error
  for a specific input to a specific benchmark,
  we consider that a true positive
  if the Rival ground-truth value differed from
  the computed floating-point value by more than
  \filled{16} ULPs.%
\footnote{
  We technically check for more than four ``bits of error'',
    which is subtly different than 16 ULPs for subnormals.
}
We then measure the rate of false positives and negatives,
  using the standard precision and recall metrics,
  to determine the accuracy of \name as a debugger.

\begin{figure}
  \centering
  \includegraphics[scale=0.58]{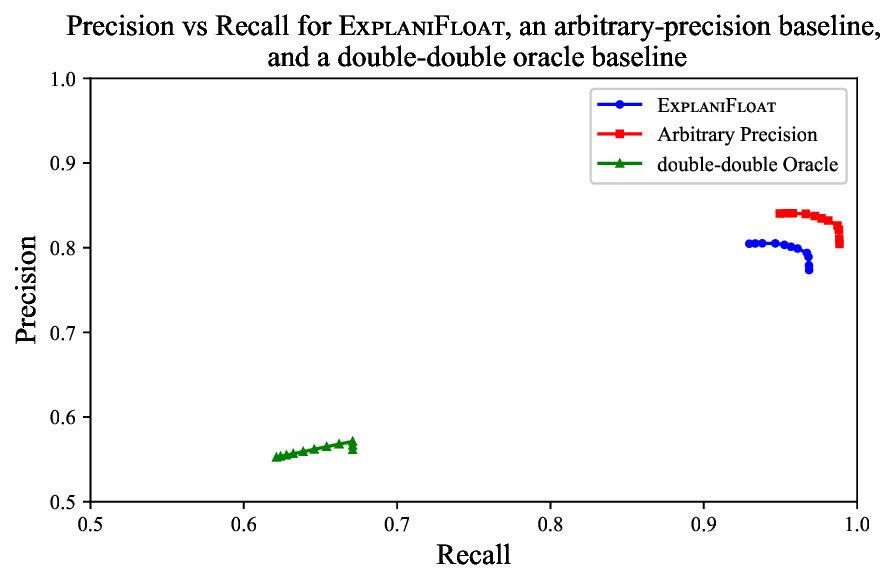}
  \caption{
    A precision vs recall graph of \name, the arbitrary-precision baseline and the double-double oracle based debugging baseline.
    We run them on thresholds starting from 4 to 4\thinspace096 doubling each time.
    Note that the precision and recall do not change drastically with change in the threshold.
  }
  \label{fig:prerec}
\end{figure}

The precision-recall plot in \Cref{fig:prerec}
  plots \name's results with a blue line.
Each point along this line shows
  the precision (vertical) and recall (horizontal)
  of \name over all inputs to all \numbench benchmarks;
  up and to the right is better.
Different points along the line use
  different condition number thresholds
  from $4$ to $4\thinspace096$.
The exact precision and recall vary by threshold,
  but for a threshold of $64$, \name has
  a precision of approximately \filled{80.0\%}
  and a recall of approximately \filled{96.1\%}.
The high recall is critical in a debugging tool:
  it means no false negatives that would hide
  the true source of error.
The high but lower precision, by contrast, is less of a concern
  because a debugging tool is typically used only when
  a problem of some kind is already known to occur.
Note that different thresholds have
  similar precision and recall results;
  this means that users do not have to fine-tune the threshold
  to get good results from \name.

To concretize these results,
  we examine four specific benchmarks
  where \name either performs well
  or suffers from false positives and negatives.
On the ``Asymptote C'' benchmark,
%https://nightly.cs.washington.edu/reports/herbie/1712320797:nightly:explanations-eval:0a13acd6df/mathematics/37-AsymptoteC/graph.html
\[
\frac{x}{x + 1} - \frac{x + 1}{x - 1}
\text{ for } x = -1.3337344672928248\cdot10^{72},
\]
  \name has perfect accuracy and recall:
  it detects a cancellation issue for \filled{117} of 256 inputs,
  and in exactly those \filled{117} cases the floating-point error surpasses 16 ULPs.

In the ``HairBSDF, Mp, Lower'' benchmark
  drawn from a computer graphics textbook~\cite{mpfr},%
%https://nightly.cs.washington.edu/reports/herbie/1711195957:nightly:explanations-eval:875f2f51aa/graphics/15-HairBSDFMplower/graph.html
\[
  \exp\left(\left(\left(\left(\frac{c_i c_O}{v} - \frac{s_i 
      s_O}{v}\right) - \frac{1}{v}\right) + 0.6931\right) + \log
    \left(\frac{1}{2 v}\right)\right),
\]
  \name's handling of overflow and underflow is essential.
For specific, unusual inputs,
  the subexpression $s_i \cdot s_O$ underflows,
  but dividing by $v$ brings the result back into range.
However, the resulting value $(s_i \cdot s_O) / v$
  is then added to $(c_i \cdot c_O) / v$, a much larger number,
  so whether or not $(s_i  \cdot s_O) / v$ underflows has minimal impact
  on the overall expression's floating-point error,
  and \name suppresses the underflow explanation
  and avoids generating a false positive.
Across all inputs to this benchmark,
  underflow suppression reduces the number of false positives
  from \filled{112} to \filled{1}.

Meanwhile, ``Expression, p6'',
%https://nightly.cs.washington.edu/reports/herbie/1712320797:nightly:explanations-eval:0a13acd6df/numerics/31-Expressionp6/graph.html
\(
  \left(a + \left(b + \left(c + d\right)\right)\right) 2
\),
  has a false negative for specific a input
  $a \approx -13.58\ldots$,
  $b \approx -2.32\ldots$,
  $c \approx 3.08\ldots$,
  $d \approx 12.53\ldots$,
  the final addition (between $a$ and $b + (c + d)$)
  has a condition number of about $45$.
This causes false negatives
  at higher condition number thresholds like $64$,
  though not at lower thresholds like $32$.
\filled{68} other inputs to this benchmark have
  similarly-middling condition numbers.

\name also sometimes generates false positives.
For example, consider the ``Spherical law of cosines'' benchmark,
%https://nightly.cs.washington.edu/reports/herbie/1712320797:nightly:explanations-eval:0a13acd6df/mathematics/12-Sphericallawofcosines/graph.html
\[
\cos^{-1} \left(\sin \phi_1 \cdot \sin \phi_2 + \left(\cos \phi_1 \cdot \cos \phi_2\right) \cdot \cos \left(\lambda_1 - \lambda_2\right)\right) R
\]
For inputs where $\lambda_1$ is large but $\lambda_2$ is very small,
  \name generates a false positive error for $\cos(\lambda_1 - \lambda_2)$.
Since $\lambda_1$ is large, the condition number for $\cos$ is large,
  but $\lambda_1 - \lambda_2$ introduces
  such a tiny relative error (roughly $2\cdot10^{-182}$),
  that even amplifying it by the very large condition number
  doesn't cause much end-to-end error.
Addressing this cause of false positives
  would be an interesting direction for future work.

\section{RQ2: Comparison to Oracle Method}

\name's exact precision and recall are less indicative
  than how it compares to the oracle method.
We thus compare \name to
  a double-double oracle-method debugger
  inspired by EFTSanitizer~\cite{eftsan}.
This variant of \name evaluates the program
  using double-precision shadow memory,
  just like \name,
  but detects errors by comparing
  the standard floating-point evaluation to the oracle value.

\Cref{fig:prerec} plots this oracle method baseline
  in green, for a range of ULP error thresholds.
The oracle baseline has a much worse
  precision and recall, topping out at
  a precision of \filled{56.5\%} and a recall of \filled{65.4\%}.
At no tested threshold value
  is the oracle-method debugger competitive with \name.

A closer look at the benchmarks
  shows that the issue is as expected:
  the oracle suffers from rounding error
  that masks or hides the rounding error
  in the floating-point evaluation.
Consider the ``2sqrt'' benchmark,
  drawn from a numerical method textbook~\cite{book87-nmse},
\[
\sqrt{x + 1} - \sqrt{x},
  \text{ for }x = 10^{100}.
\]
Here the oracle-method baseline
  computes the same result for both $\sqrt{x+1}$ and $\sqrt{x}$,
  resulting in a final oracle value of $0$.
The floating-point computation also computes $0$,
  meaning no error is raised;
  across all inputs to this benchmark,
  the oracle-method baseline has a recall of \filled{11.3\%}.
\name, on the other hand,
  detects a very large condition number in this case
  and achieves a perfect \filled{100.0\%} recall
  across all inputs to this benchmark.

The oracle baseline also handles overflow and underflow poorly.
Consider the benchmark ``cos2'' from the same source,
\[
(1 - \cos(x)) / (x\times x), \text{ for } x = 10^{200}
\]
The double-double oracle warns about overflow
  in $x \times x$.
However, both the correct and computed floating-point results are 0,
  meaning this input actually has no error.
\name correctly handles this case
  by computing the logarithm of the output value
  as outside the representable range.
Since the value is never brought back \emph{into} range,
  it does not produce a warning.

\subsection{RQ3: Equal performance to arbitrary precision}

Finally, we aim to show that using double-double values
  provides enough precision for accurate condition number computation.
We thus modify \name to use
  Rival's correctly-rounded arbitrary-precision baseline
  for all intermediate values,
  but still produce errors using condition numbers
  and overflow renormalization.
Because the intermediate values are computed exactly,
  this variation allows us to evaluate whether
  the use of double-double shadow values
  introduces additional false positives and negatives.

\Cref{fig:prerec} plots this baseline in red.
On average, this baseline
  achieves a precision of \filled{83.2\%}
  and a recall of \filled{98.1\%}.
These are only slightly higher (by \filled{3.2\%} and \filled{2.0\%})
  than \name,
  showing that the precision of \name's shadow values
  do not significantly affect its results.
%4.27 faster
Despite the largely-similar predictive accuracy,
  DSL is significantly faster than the alternative baseline,
  taking \filled{3.7 seconds} in \name
  versus \filled{15.7 seconds} with the alternative baseline,
  a speedup of $4.24\times$.
Since \name was not engineered for maximum performance,
  we take this speedup number with a grain of salt,
  but we do expect \name's performance advantage
  to be substantial,
  since it avoids allocation and arbitrary-precision computations,
  and to grow even larger for larger programs.

Comparing \name to the perfect-oracle baseline,
  we find that \name's shadow values introduce
  three core limitations:
  cancellation in large sums;
  residual underflow;
  and aliasing between errors
  in transcendental functions.
Cancellation in large sums refers to cases
  where three or more values are summed together
  and multiple cancellations occur,
  like in the ``exp2'' benchmark,
\[
  (e^x - 2) + e^{-x}
\]
For $x$ close to zero, such as $x = 10^{-200}$,
  the exponential terms in this expression
  implicitly act like the sum $1 \pm x + x^2/2$,
  meaning that this expression
  effectively adds seven terms,
  of which 5 (1, 1, and $-2$; $x$ and $-x$) cancel.
DSL is not effective on this input
  because double-double evaluation of $e^x$
  retains only the $1$ and $\pm x$ term.
In other words, DSL only stores the values that cancel,
  so its final evaluation of the expression is $0$.
While \name does detect a high condition number,
  it is not able to determine whether underflow occurs,
  which causes false positives.

Residual underflow refers to cases
  where the residual term in a double-double value
  cannot be represented in floating-point,
  but the primary term can.
For example, consider the
  the ``sintan'' benchmark:
\[(x - \sin(x)) / (x - \tan(x))\]
For $x$ very close to zero,
  such as $x = 10^{-200}$,
  $\sin(x)$ and $\tan(x)$
  evaluate to $x$
  with a residual value of approximately
  $10^{-600}$.
However, this residual value in fact underflows,
  meaning that in effect \name performs
  only a double-precision evaluation of the benchmark.
Here, \name does correctly warn due to the high condition number,
  but does not also produce a renormalization error
  because it cannot compute the exponent of the subtraction.
This issue was also noted in EFTSanitizer~\cite{eftsan},
  but the issue was rare in that paper's evaluation
  on mostly-linear-algebra workloads.
In our larger and more diverse benchmark suite,
  it does cause false negatives.

Aliasing refers to cases where one operation's rounding error
  is cancelled by another operation's rounding error.
For example, consider the ``logs'' benchmark:
\[((n + 1)\log(n + 1) - n \log n) + 1\]
For large $n$, like $n = 10^{200}$,
  $\log(n + 1)$ and $\log(n)$
  are very close and have nearly-identical rounding error.
In this case, $(n + 1) \log(n + 1)$ and $n \log n$
  have the same shadow value
  and so their difference evaluates to exactly 0.
\name raises a condition number error, but the later addition to $1$
  means that the error is (incorrectly) suppressed,
  leading to a false negative.
In reality, the error is approximately $\log n$,
  much larger than $1$.
EFTSanitizer likely did not have this issue
  due to its limited support for transcendental functions.

All that said, \name's precision and recall are
  very similar to the arbitrary-precision baseline,
  showing that the performance benefits of \name's shadow values
  come with very few downsides in precision or recall.

%\subfile{relwork}

\section{Conclusion}
\label{sec:conclusion}

\name combines recent advances
  in numeric debugging and static analysis tools
  to create an accurate yet performant numerical debugger.
It uses condition numbers instead of oracles to detect rounding error
  and uses a novel oracle for detecting over- and underflows.
The result has exceptional precision
  (\filled{80.0\%}) and recall (\filled{96.1\%}),
  beating both double-double oracle and arbitrary-precision approaches.

\enlargethispage{-.4in}
\bibliographystyle{IEEEtranN}
\bibliography{references}

\end{document}